\documentclass[11pt]{article}
\usepackage{amssymb,amsfonts,amsmath,latexsym,epsf,tikz,url}

\newtheorem{theorem}{Theorem}[section]
\newtheorem{proposition}[theorem]{Proposition}

\newtheorem{corollary}[theorem]{Corollary}

\newtheorem{remark}[theorem]{Remark}
\newtheorem{definition}[theorem]{Definition}

\newcommand{\proof}{\noindent{\bf Proof.\ }}
\newcommand{\qed}{\hfill $\square$\medskip}

\begin{document}

\title{The cost number and the determining number of a graph}

\author{
Saeid Alikhani$^{}$\footnote{Corresponding author}  \and Samaneh Soltani 
}

\date{\today}

\maketitle

\begin{center}
Department of Mathematics, Yazd University, 89195-741, Yazd, Iran\\
{\tt alikhani@yazd.ac.ir, s.soltani1979@gmail.com}
\end{center}

\begin{abstract}
 
The distinguishing number $D(G)$ of a graph $G$ is the least integer $d$ such that
$G$ has an vertex labeling with $d$ labels that is preserved only by a trivial automorphism. The minimum size of a label class in such a labeling of $G$ with $D(G) = d$
is called the cost of $d$-distinguishing $G$ and is denoted  by $\rho_d(G)$.
A set of vertices $S\subseteq V(G)$ is a determining set for $G$ if every automorphism of $G$ is uniquely determined by
its action on $S$. The determining number of $G$, ${\rm Det}(G)$,  is the minimum cardinality of determining sets of $G$. 
 In this paper we
obtain some general upper and lower bounds for $\rho_d(G)$ based on ${\rm Det}(G)$. Finally, we compute the
cost and the determining number for the  friendship graphs and corona product of two graphs.
\end{abstract}


\medskip
\noindent{\bf AMS Subj.\ Class.:} 05C15, 05C25

\section{Introduction}

Let $G = (V,E)$ be a simple graph with $n$ vertices. We use the standard graph notation
(\cite{Sandi}). The set of all automorphisms of $G$, with the operation of composition of permutations, is a permutation group on $V$ and is denoted by ${\rm Aut}(G)$. A labeling of $G$,
 $\phi : V \rightarrow \{1, 2, \ldots , r\}$, is \textit{$r$-distinguishing}, if no non-trivial automorphism of $G$ preserves
all of the vertex labels. In other words, $\phi$ is $r$-distinguishing if for every non-trivial
$\sigma \in {\rm Aut}(G)$, there exists $x$ in $V$ such that $\phi(x) \neq  \phi(\sigma x)$. The  \textit{distinguishing number} of a graph $G$ has been defined in \cite{Albert} and is the minimum
number $r$ such that $G$ has a labeling that is $r$-distinguishing. We will also need to
know what it means for a subset of vertices to be $d$-distinguishable. For  $W \subseteq V (G)$, a labeling $f : W \rightarrow \{1, \ldots ,  d\}$ is called \textit{$d$-distinguishing} if whenever an automorphism
fixes $W$ setwise and preserves the label classes of $W$ then it fixes $W$ pointwise. Note that
though such an automorphism fixes $W$ pointwise, it is not necessarily trivial; it may
permute vertices in the complement of $W$. A set $W$ is called \textit{$d$-distinguishable} if it has a
$d$-distinguishing labeling. By definition, $W$ is $1$-distinguishable if every automorphism
that preserves $W$ fixes it pointwise. The introduction of the distinguishing number was
a great success; by now about one hundred papers were written motivated by this seminal paper! The core of the research has been done on the invariant $D$ itself, either on
finite \cite{Chan, immel-2017, Kim} or infinite graphs \cite{estaji-2017, lehner-2016, smith-2014}; see also the references therein.

In 2007 Wilfried Imrich posed the following question \cite{wilfImr}: ``What is the minimum
number of vertices in a label class of a $2$-distinguishing labeling for the hypercube $Q_n$?" To aid in addressing this question,
Boutin \cite{Boutin2008} called a label class in a $2$-distinguishing labeling of $G$ a \textit{distinguishing class}. She called the minimum size of such a class in $G$ the \textit{cost of $2$-distinguishing $G$} and
denoted it by $\rho(G)$. Boutin also showed that $\lceil {\rm log}_2n\rceil - 1 \leq  \rho (Q_n) \leq \lceil {\rm log}_2n\rceil + 1$.
She used the \textit{ determining set} \cite{Boutin2006}, a set of vertices whose pointwise stabilizer is trivial.
In other words, a subset $S$ of the vertices of a graph $G$ is called a determining set if
whenever $g,h \in {\rm Aut}(G)$ agree on the vertices of $S$, they agree on all vertices of $G$. That
is, $S$ is a determining set if whenever $g$ and $h$ are automorphisms with the property that
$g(s) = h(s)$ for all $s\in S$,  then $g = h$. 
Albertson and Boutin proved the following theorem in \cite{Boutin2006}. 

\begin{theorem}{\rm \cite{Boutin2006}}\label{1.1} 
	 A graph is $d$-distinguishable if and only if it has a determining set
	that is $(d-1)$-distinguishable.
\end{theorem}
 In particular, the complement of
such a determining set is a label class in a $d$-distinguishing labeling of $G$. Thus, a graph
is $2$-distinguishable if and only if it has a determining set for which any automorphism
that fixes it setwise must also fix it pointwise. In such a case, the determining set and
its complement provide the two necessary label classes for a $2$-distinguishing labeling.
Thus, in particular, the cost of $2$-distinguishing a graph $G$ is bounded below by the size
of a smallest determining set, denoted ${\rm Det}(G)$.

\medskip

In this paper the cost of $2$-distinguishing is extended to the cost of $d$-distinguishing.
This paper is organized as follows. Definitions and facts about the cost of $d$-distinguishing,
is given in Section 2. Also by finding the cost number and the determining number of the friendship graph in Section 2, we show that for any positive integer $m$, there exists a graph $G$ with $D(G) = d$ such
that $|{\rm Det}(G)-\rho_d(G)| = m$.
 The cost of $d$-distinguishing corona product of two graphs are given in
Section 3.

\section{The cost of $d$-distinguishing graphs}
We start with the following definition: 
\begin{definition}\label{2.1} Let $G$ be a graph with the distinguishing number $D(G) = d$. The
minimum size of a label class in any $d$-distinguishing labeling of $G$, is called the cost of
$d$-distinguishing of $G$ and denoted it by $\rho(G)$, or by $\rho_d(G)$ if we wish to stress that the
distinguishing number of $G$ is $d$.
\end{definition}

It can be easily seen that the cost of $n$-distinguishing of complete graph $K_n$ and
complete bipartite graph $K_{n,m}$ $(m < n)$ and $K_{n-1,n-1}$ is $1$. The following result is an 
immediate consequence of Definition \ref{2.1}.
\begin{proposition}\label{2.2} 
Let $G$ be a graph of order $n$ and the distinguishing number $D(G) = d$.
Then
\begin{enumerate}
\item[(i)] The cost of $d$-distinguishing graph $G$ is $\rho_d(G) \leq \frac{n}{d}$.
\item[(ii)] $d = 1$ if and only if $\rho_d(G) = n$.
\end{enumerate}
\end{proposition}

\begin{proposition}\label{2.3}  
If $G$ is a graph of order $n$ with  the distinguishing number $D(G) =d \geq 2$, then $\rho_d(G) \leq  \frac{n}{2}$. In particular, if $\rho_d(G) =  \frac{n}{2}$, then $d = 2$
\end{proposition}
\proof The first part follows directly from Proposition \ref{2.2}. For the second part, since
$\rho_d(G) =  \frac{n}{2}$, so the size of the remaining label classes in $d$-distinguishing labeling is at
least $\frac{n}{2}$, thus we have exactly two distinguishing classes, and hence $D(G) = d = 2$.\qed

The converse of Proposition \ref{2.3} is not true, for instance see the path graphs $P_n$
with $D(P_n) = 2$ and $\rho_2(P_n) = 1$ where $n \geq 3$.
Here we obtain some bounds for the cost of $d$-distinguishing graphs using its determining number.

\begin{proposition}\label{2.4}  
Let $G$ be a graph with the distinguishing number $D(G) = d$. If $\psi$  is
a $d$-distinguishing labeling of $G$ with distinguishing classes of sizes $t_1\leq \cdots \leq  t_d$ such
that $t_1 = \rho_d(G)$, then
$${\rm Det}(G) \leq \rho_d(G) + t_2 + \cdots + t_{d-1}.$$
\end{proposition}
\proof Since the union of all distinguishing classes of sizes $t_1, \ldots , t_{d-1}$ is a determining
set of $G$ and $t_1 = \rho_d(G)$, so we have the result. \qed

The upper bound of Proposition \ref{2.4} is sharp for complete graphs and star graphs.
\begin{proposition}\label{2.5}  
Let $G$ be a graph of order $n$ and the distinguishing number $D(G) = d$.
Then $\rho_d(G) \leq n - {\rm Det}(G)$.
\end{proposition}
\proof The distinguishing number of a determining set of size ${\rm Det}(G)$ is at most $d-1$,
by Theorem \ref{1.1}. Since  the complement
of such a determining set is a label class in $d$-distinguishing labeling of $G$, so we have the result. \qed

\begin{corollary} \label{2.6}
 Let $G$ be a graph of order $n$ and  $D(G) = d$.
If the distinguishing number of a determining set of $G$ of size ${\rm Det}(G)$, say $A$, is $d- 1$,
then
$$\rho_d(G) \leq {\rm min}\{n - {\rm Det}(G), \rho_{d-1}(G[A])\},$$
where $G[A]$ is the induced subgraph of $G$ generated by vertices in $A$.
\end{corollary}
\proof Set ${\rm Det}(G) = t$ and let $A = \{v_1, \ldots , v_t\}$ be a determining set of $G$ with the
distinguishing number $d - 1$. If we label the vertices of $G[A]$ with labels $1, \ldots , d - 1$
distinguishingly, and label all vertices $v_{t+1},\ldots , v_n$ with new label $d$, then it can be seen
that we have a distinguishing labeling of $G$ with $d$ labels. Since the minimum size of
distinguishing classes of this labeling is ${\rm min}\{n-t, \rho_{d-1}(G[A])\}$, so we have the result.\qed

By Proposition \ref{2.5} and the fact that ${\rm Det}(G) \leq  \rho_2(G)$, we can prove the following
result.
\begin{corollary} \label{2.7}
 Let $G$ be a graph of order $n$ and the distinguishing number $D(G) = d$.
\begin{enumerate}
\item[(i)] If ${\rm Det}(G) \leq \rho_d(G)$, then ${\rm Det}(G) \leq \frac{n}{2}$.
\item[(ii)] If $d = 2$, then ${\rm Det}(G) \leq \frac{n}{2}$.
\end{enumerate}
\end{corollary}

\begin{figure}[ht]
	\hspace{1cm}
	\begin{minipage}{6.3cm}
		\includegraphics[width=\textwidth]{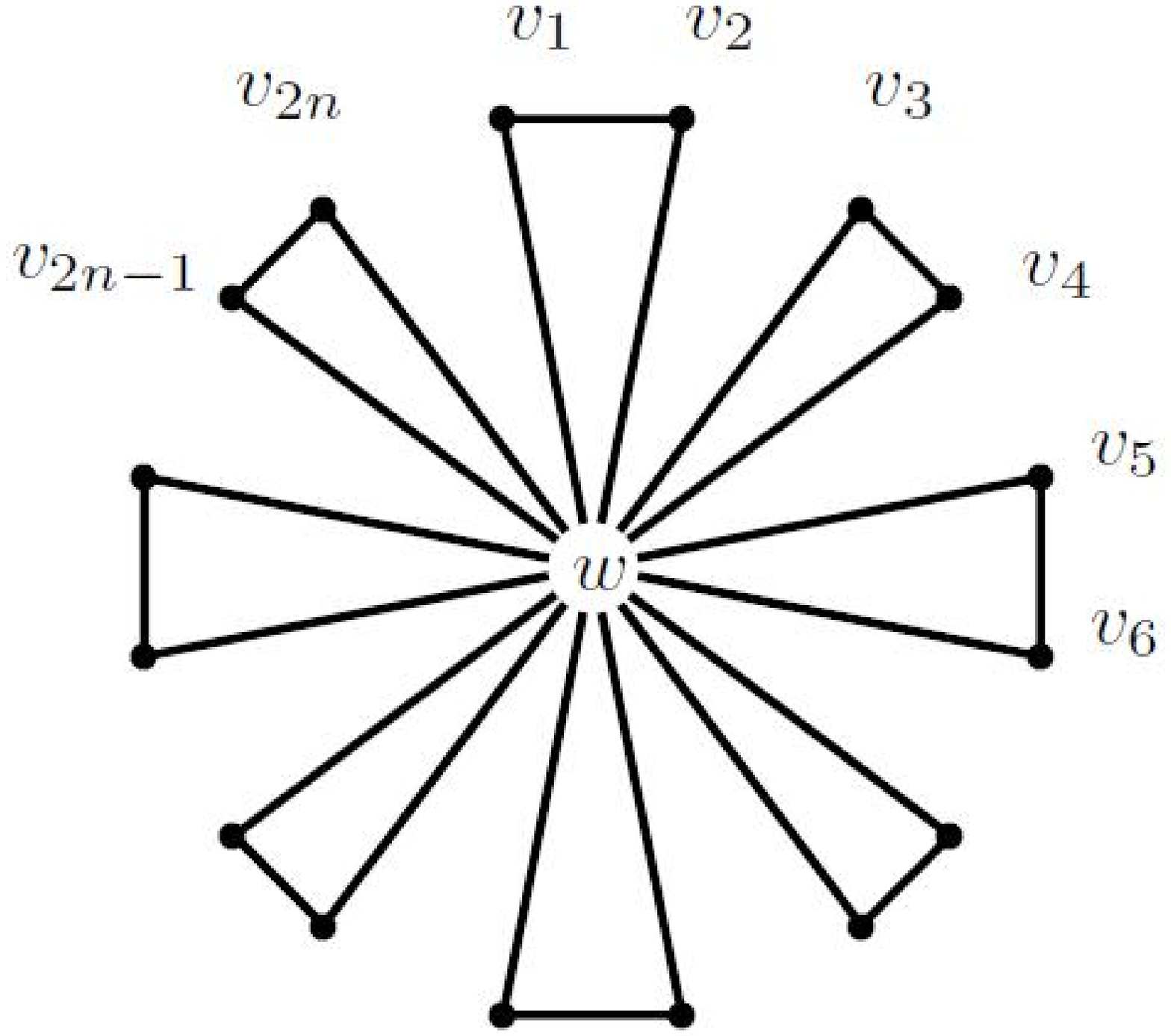}
	\end{minipage}
	\begin{minipage}{6.1cm}
		\includegraphics[width=\textwidth]{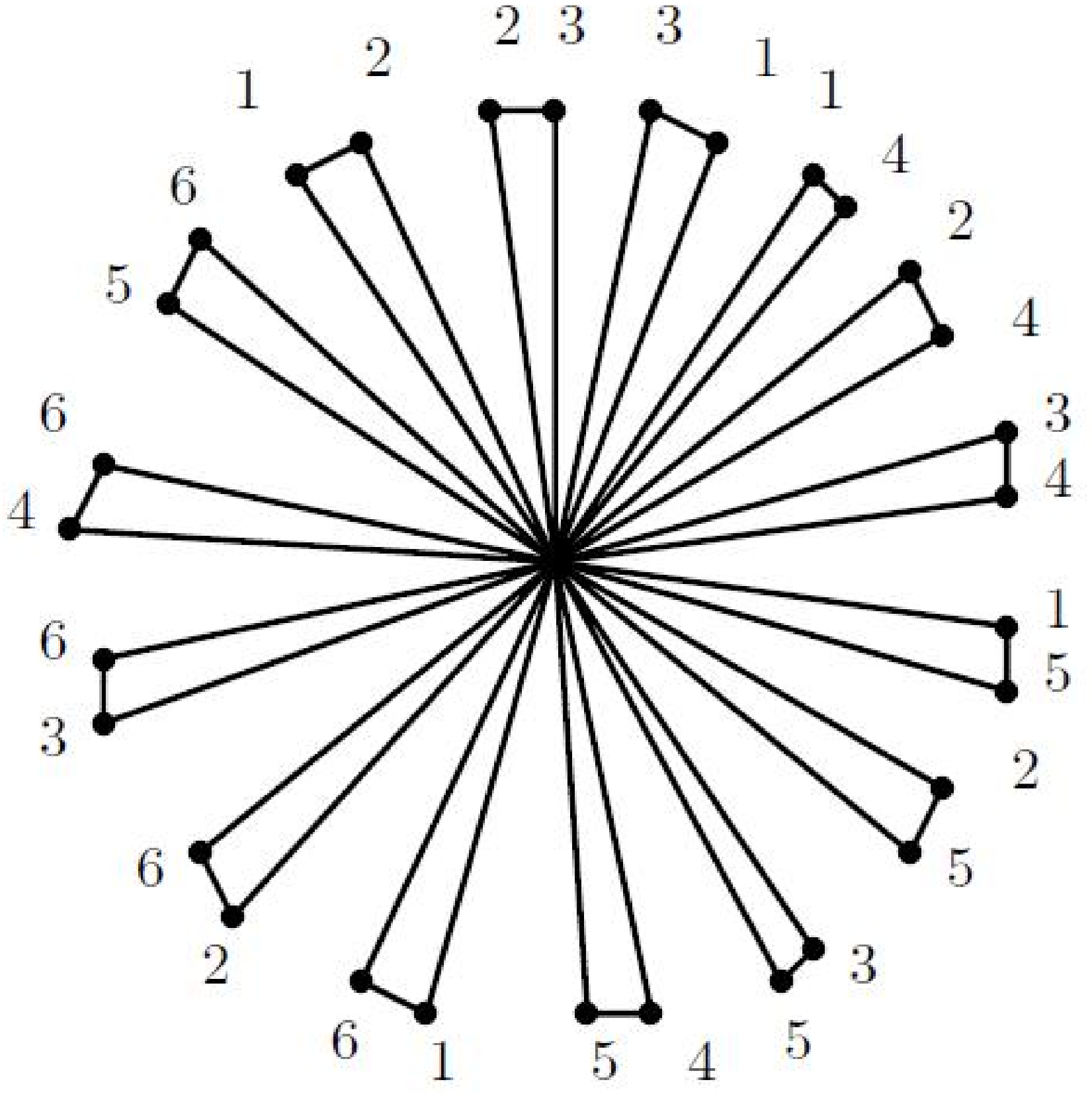}
	\end{minipage}
	\caption{\label{friend} Friendship graph $F_n$ and the vertex labeling of $F_{15}$, respectively.}
\end{figure}

We shall show that for any positive integer $m$, there exists a graph $G$ with $D(G) = d$ such
that $|{\rm Det}(G)-\rho_d(G)| = m$. To do this 
 we consider the friendship graphs and compute their cost and determining
number. The friendship graph $F_n$ $(n \geq 2)$ can be constructed by joining $n$ copies of
the cycle graph $C_3$ with a common vertex (see Figure \ref{friend}). The authors obtained the
distinguishing number of friendship graphs as follows:

\begin{theorem}{\rm \cite{Alikhani}}\label{3.1}
The distinguishing number of the friendship graph $F_n$ $(n \geq 2)$ is
$$D(F_n) = \lceil \dfrac{1+\sqrt{8n+1}}{2}\rceil.$$
\end{theorem}

\begin{remark}\label{3.2}
 Let $k_j = {\rm min}\{i : ~D(F_i) = j\}$ for any $j \geq 3$. Then by Theorem \ref{3.1} we
have:
\begin{enumerate}
\item  For any $j \geq 3$, $k_j = \lfloor\frac{j^2-3j+2}{2}\rfloor +1$.
\item  For all $i$, $0 \leq i \leq j - 2$, $D(F_{k_j+i}) = j$ and $D(F_{k_j+j-1}) = j + 1$.
\end{enumerate}
\end{remark}
\begin{theorem}\label{3.3}
  Let $j \geq 3$ and $k_j = {\rm min}\{i : D(F_i) = j\}$. Then $\rho_j(F_{k_j+i}) = i + 1$ where
$0 \leq i \leq j - 2$.
\end{theorem}
\proof In any $j$-distinguishing labeling of $F_{k_j+i}$ with labels $\{1, \ldots , j\}$, each of the 2-sets
consisting of vertex of degree two and its neighbor of degree two must have a different
2-subset of labels $\{1,\ldots , j\}$. Since $k_j = {\rm min}\{i : D(F_i) = j\}$, so the all 2-subsets of $\{1,\ldots j\}$ have been used for any distinguishing labeling of $F_{k_j-1}$. Thus without loss of
generality, we can assume that the number of label $p$ which is used for labeling of vertex
set of $F_{k_j-1}$, say $n_p(F_{k_j-1})$, is $n_p(F_{k_j-1}) = j-2$ for $2 \leq p \leq j-1$ and $n_1(F_{k_j-1}) = j-1$
(the central vertex w is labeled with label 1). If we assign the 2-sets $\{v_{2q-1}, v_{2q}\}$, where
$k_j \leq q \leq k_j + i$, the 2-subsets $\{i + 1, j\}$ of labels, then we obtain a distinguishing
labeling for $F_{k_j+i}$ with labels $1, \ldots , j$ such that
$$n_j(F_{k_j+i}) = {\rm min}\{n_1(F_{k_j+i}),\ldots , n_j(F_{k_j+i})\} = i + 1.$$
 Thus $\rho_j(F_{k_j+i}) \leq i + 1$. On the other hand, we have $n_p(F_{k_j-1}) \geq j - 2$, for any
$2 \leq p \leq j - 1$, so since $n_p(F_{k_j+i}) \geq  n_p(F_{k_j-1}) \geq j - 2$ and $0 \leq i \leq j - 2$, we
have $n_j(F_{k_j+i}) = {\rm min}\{n_1(F_{k_j+i}), \ldots , n_j(F_{k_j+i})\}$. Now since the label $j$ have been used
only for vertices $v_q$, where $2k_j - 1 \leq q \leq 2k_j + 2i$, and since the 2-subsets of labels
related to the 2-sets $\{v_{2q-1}, v_{2q}\}$ and $\{v_{2q'-1}, v_{2q'}\}$ must be different for any $q, q' \in \{k_j , k_j+1,\ldots , k_j+i\}$ where $q \neq q'$, so $n_j(F_{k_j+i}) = i+1$, and therefore $\rho_j(F_{k_j+i}) = i+1$.\qed

\begin{theorem}\label{3.4}
  For any $n \geq 2$, ${\rm Det}(F_n) = n$.
\end{theorem}
\proof Let the vertices of $F_n$ be as shown in Figure \ref{friend}. It can be easily seen that the
set $A = \{v_1, v_3,\ldots ,v_{2n-1}\}$ is a determining set for $F_n$. On the other hand, if $B$ is
a determining set of $F_n$ with $|B| \leq n - 1$, then there exists $i \in \{1, \ldots , n\}$ such that
$v_{2i-1}, v_{2i} \notin B$. Hence there exists the nonidentity automorphism $f$ of $F_n$ with $f(x) = x$
for all $x \in B$, $f(v_{2i-1}) = v_{2i}$ and $f(v_{2i}) = v_{2i-1}$, which is a contradiction to that $B$ is a
determining set. Therefore ${\rm Det}(F_n) = n$. \qed

Now we end this section by the following theorem: 
\begin{theorem}\label{2.8} 
	For any positive integer $m$, there exists a graph $G$ with $D(G) = d$ such
	that $|{\rm Det}(G) - \rho_d(G)| = m$.
\end{theorem}

\proof 
 By Theorems \ref{3.3} and \ref{3.4}, it can be concluded that for every
positive integer $m$, there exists some suitable $n$ such that the friendship graph $F_n$
satisfies $|{\rm Det}(F_n) - \rho_d(F_n)| = m$. \qed

\section{The cost and determining number of corona product}
In this section, we shall study the cost number  and the determining number of corona product of
graphs. The \textit{corona product}  $G \circ H$ of two graphs $G$ and $H$ is defined as the graph
obtained by taking one copy of $G$ and $|V (G)|$ copies of $H$ and joining the $i$-th vertex of
$G$ to every vertex in the $i$-th copy of $H$. The distinguishing number of corona product
of graphs have been studied by the authors in \cite{Alikhani}. Before presenting our results, we
explain the relationship between the automorphism group of the graph $G \circ H$ with
the automorphism groups of two connected graphs $G$ and $H$ such that $G \neq K_1$. Note
that there is no vertex in the copies of $H$ which has the same degree as a vertex in
$G$. Because if there exists a vertex $w$ in one of the copies of $H$ and a vertex $v$ in
$G$ such that ${\rm deg}_{G\circ H} (v) = {\rm deg}_{G\circ H} (w)$, then ${\rm deg}_G(v) + |V (H)| = {\rm deg}_H(w) + 1$. So we
have  ${\rm deg}_H(w) + 1> |V (H)|$, which is a contradiction. Let the vertex set of $G$ be
$\{v_1, \ldots , v_{|V (G)|}\}$ and the vertex set of $i$-th copy of $H$, $H_i$, be $\{w_{i1}, \ldots , w_{i|V (H)|}\}$. Since
there is no vertex in copies of $H$ which has the same degree as a vertex in $G$, for
every $f \in  {\rm Aut}(G  \circ H)$, we have $f|_H \in {\rm Aut}(H)$ and $f|_G \in {\rm Aut}(G)$. In addition, for
$i, j \in \{1, \ldots , |V (G)|\}$ we have
$$f(v_i) = v_j \Longleftrightarrow f(H_i) = H_j. $$
Conversely, let $\varphi \in {\rm Aut}(G)$ and $\phi \in {\rm Aut}(H)$ such that $\varphi(v_i) = v_{j_i}$, where $i, j_i  \in \{1, \ldots , |V (G)|\}$. Now we define the following automorphism $h$ of $G \circ H$:
\begin{equation*}
h:\left\{
\begin{array}{ll}
v_i \mapsto \varphi(v_i) = v_{j_i} &i, j_i \in \{1,\ldots , |V (G)|\},\\
w_{ik} \mapsto \phi(w_{j_ik}) &k \in \{1,\ldots , |V (H)|\}.
\end{array}\right.
\end{equation*}

We start with the determining number of corona product of two graphs.
\begin{theorem}\label{4.1}
 Let $G$ and $H$ be two connected graphs of orders $n,m \geq 2$, respectively.
Then
$${\rm Det}(G \circ H) = {\rm Det}(G) + n{\rm Det}(H).$$
\end{theorem}
\proof We denote the vertices of $G$ in $G \circ H$ by $v_1, \ldots , v_n$, and vertices of $H$ corresponding to the vertex $v_i$ by $w_{i1}, \ldots , w_{im}$. Let ${\rm Det}(G) = k$ and ${\rm Det}(H) = k'$. We
suppose that the sets $\{v_1,\ldots , v_k\}$ and $\{w_1,\ldots , w_{k'}\}$ are the determining sets of $G$
and $H$, respectively, then the set $\{v_1,\ldots, v_k\}\cup  (\bigcup_{i=1}^n \{w_{i1}, \ldots , w_{i{k'}}\})$ is a determining
set of $G \circ H$, and hence ${\rm Det}(G \circ H) \leq {\rm Det}(G) + n{\rm Det}(H)$. On the other hand if
${\rm Det}(G \circ H) < {\rm Det}(G) + n{\rm Det}(H)$, then there exists a determining set $Z$ for $G\circ H$ with
$|Z| = {\rm Det}(G \circ H)$ such that $|Z \cap V (H_i)| < k'$ or $|Z \cap V (G)| < k$ for some $1 \leq i \leq n$,
where $H_i$ is the isomorphic copy of $H$ corresponding to the vertex $v_i$ in $G \circ H$. We
consider the two following cases:
\begin{enumerate}
\item[Case 1)] Let $Z \cap V (H_i) = \{w_{ij_1}, \ldots , w_{ij_t}\}$ where $t < k'$ for some $i$, $1 \leq i \leq n$. Since
$t < k'$, it can be concluded that there exists a nonidentity automorphism $f$ of
$H$ such that $f(w_{ij_1}) = w_{ij_1},\ldots , f(w_{ij_t}) = w_{ij_t}$. We extend $f$ to a nonidentity
automorphism $\overline{f}$ of $G \circ H$ with
\begin{equation*}
\overline{f}(x) =\left\{
\begin{array}{ll}
x & \text{if}~ x \in V (G),\\
f(x) & \text{if} ~x \in V (H_i),\\
x &\text{if}~ x \in V (H_{i'}), i' \neq i.
\end{array}\right.
\end{equation*}

In this case, $\overline{f}$ is a nonidentity automorphism of $G\circ H$ and it fixes the determining
set $Z$, pointwise, which is a contradiction.
\item[Case 2)] Let $Z \cap V (G) = \{v_{j_1} ,\ldots , v_{j_t}\}$ where $t < k$. Since $t < k$, so there exists a
nonidentity automorphism $f$ of $G$ such that $f(v_{j_1}) = v_{j_1},\ldots , f(v_{j_t}) = v_{j_t}$. We
extend $f$ to a nonidentity automorphism $\overline{f}$ of $G \circ  H$ with
\begin{equation*}
\overline{f}(x) =\left\{
\begin{array}{ll}
f(x) & \text{if}~ x \in V (G),\\
x &\text{if}~ x \in V (H_{i}), i=1,\ldots , n.
\end{array}\right.
\end{equation*}

In this case, $\overline{f}$ is a nonidentity automorphism of $G\circ H$ and it fixes the determining
set $Z$, pointwise, which is a contradiction. \qed
\end{enumerate}

\begin{theorem}\label{4.2}
 If $G$ is a connected graph of order $n \geq 2$, then ${\rm Det}(G\circ K_1) = {\rm Det}(G)$.
\end{theorem}
\proof It is clear that each determining set of $G$ is a determining set of $G \circ  K_1$,
and so ${\rm Det}(G \circ  K_1) \leq {\rm Det}(G)$. Set ${\rm Det}(G) = k$, $V (G) = \{v_1,\ldots , v_n\}$, and denote the vertex of $K_1$ adjacent to the vertex $v_i$, by $w_i$. Assume by contrary that
$t ={\rm  Det}(G \circ  K_1) < k$. Then, there exists a determining set $Z$ of $G \circ  K_1$ such that $Z =\{v_1,\ldots , v_{t_1} ,w_{j_1} ,\ldots ,w_{j_{t-t1}}\}$ where $t_1 \leq t < k$. We show that $\{v_1, \ldots ,v_{t_1} , v_{j_1} ,\ldots , v_{j_{t-t1}}\}$ 
is a determining set of $G$ with less than $k$ elements, which is a contradiction. Before
it, we note that since $Z$ is a determining set, so $\{1,\ldots , t_1\} \cap \{j_1,\ldots , j_{t-t_1}\}= \emptyset$, since
otherwise if $j_x \in \{1,\ldots , t_1\} \cap \{j_1,\ldots , j_{t-t_1}\}$, then $Z' = Z -\{ w_{j_x}\}$ is a determining set
of $G\circ K_1$ with $|Z'| < |Z|$, which is a contradiction. If $f$ is a nonidentity automorphism
of $G$ with $f(v_i) := v_{\sigma(i)}$, where $\sigma$ is a nonidentity permutation of $1,\ldots , n$, fixing the
vertices of $\{v_1,\ldots ,v_{t_1} , v_{j_1},\ldots , v_{j_{t-t_1}}\}$, pointwise, then we can extend $f$ to the nonidentity automorphism $\overline{f}$ of $G\circ K_1$ with definition $\overline{f}(v_i) := v_{\sigma(i)}$ and $\overline{f}(w_i) = w_{\sigma(i)}$ for every
$1 \leq i \leq n$. Thus $\overline{f}$ fixes the vertices of $Z$ pointwise, which is a contradiction. Thus the
vertices of $\{v_1,\ldots ,v_{t_1} , v_{j_1},\ldots , v_{j_{t-t_1}}\}$
 is a determining set of $G$. \qed

\begin{theorem}\label{4.3}
 Let $G$ and $H$ be two connected graphs of orders $n,m \geq 2$, respectively,
with $D(G) = k$ and $D(H) = k'$. If $k'' = {\rm max}\{k, k'\}$ and $D(G \circ H) = k''$, then
$$\rho_{k''}(G \circ H) \leq \rho_k(G) + n\rho_{k'}(H).$$
\end{theorem}
\proof We present a distinguishing labeling for $G \circ H$ with $k''$ labels such that the minimum size of a distinguishing class in this $k''$-distinguishing labeling is $\rho_k(G)+n\rho_{k'}(H)$.
For this purpose, we label the vertices of $G$ distinguishingly with $k$ labels $1,\ldots , k$ such
that the distinguishing class 1 has the minimum size among others. Then we label
each of copies of $H$ distinguishingly with $k'$ labels $1,\ldots , k'$ such that the distinguishing
class 1 has the minimum size among the remaining distinguishing classes of $H$. This
labeling of $G \circ H$ is a $k''$-distinguishing labeling. In fact, if $f$ is an automorphism of
$G \circ H$ preserving the labeling, then since the restriction of $f$ to $G$ and each copy of $H$
is an utomorphism of $G$ and $H$, respectively, and since the vertices of $G$ and each copy
of $H$ have been labeled distinguishingly, so these restrictions are identity, and hence $f$
is the identity automorphism of $G \circ H$.
Since the distinguishing class 1 has the minimum size $\rho_k(G) + n\rho_{k'}(H)$ among the
remaining distinguishing classes of $G \circ H$, so the result follows. \qed

\end{document}